\documentclass[10pt,leqno,oneside]{article}
\usepackage{amsmath, amssymb}
\usepackage{amsfonts}
\newcommand{\field}[1]{\mathbb{#1}}

\newcommand{\n}{\noindent}
\newcommand{\be}{\begin{equation}}
\newcommand{\ee}{\end{equation}}
\newcommand{\ben}{\begin{displaymath}}
\newcommand{\Su}{{\mathcal{S}}}
\newcommand{\een}{\end{displaymath}}
\newcommand{\ep}{\hspace{\stretch{1}}$\Box$}

\newcommand{\vs}{\vspace{0.2cm}}

\newcommand{\Reg}{{\mathcal{R}}_{\ell}}
\newcommand{\BReg}{\bar{{\mathcal{R}}}_{\bar{\ell}}}

\usepackage{epsfig}
\usepackage{latexsym}
\usepackage{mathrsfs}

\newcommand{\vru}{\bar{\nu}}
\newcommand{\vrd}{\underline{\nu}}

\newtheorem{Definition}{Definition}
\newtheorem{Proposition}{Proposition}
\newtheorem{Theorem}{Theorem}

\newtheorem{Lemma}{Lemma}

\newtheorem{Corollary}{Corollary}

\begin{document}

\begin{center}
\Large \bf Geometric relations of stable minimal surfaces and applications.
\end{center}

\begin{center}
{\large Martin Reiris.}\footnote{e-mail: martin@aei.mpg.de}\\

\vs
\textsc{Albert Einstein Institute - Max Planck.}\\

\vspace{0.7cm}
\begin{minipage}[c]{10.5cm}
\begin{center}
{\bf\Large Abstract}
\end{center}
We establish some a priori geometric relations on stable minimal surfaces lying inside three-manifolds with scalar curvature uniformly bounded below. The
relations are based on a slight generalization of a formula due to Castillon. We apply it to prove non-local rigidity results in the particular sense that they express how 
local isoperimetric properties in some region affect the local isoperimetric properties in any other region. We present applications to understand the notion of isoperimetric collapse
on three-manifolds with scalar curvature uniformly bounded below. 
\end{minipage}
\end{center}

\vspace{0.5cm}
\begin{center}
\begin{minipage}[c]{12.5cm}{\small
{\center\tableofcontents}}
\end{minipage}
\end{center}

\newpage
{\center \section{Introduction.}}
\vs
It has been observed that the geometry of stable minimal surfaces in three-manifolds is  rather sensitive on the scalar curvature of the ambient space. Therefore once one shows a priori geometric relations, these stable surfaces become automatically interesting as tools to explore the geometry of the scalar curvature. With this in mind, we establish a set of geometric properties, some of them of isoperimetric nature, for stable minimal surfaces inside three-manifolds with scalar curvature uniformly bounded below which depend on the curvature of the ambient space only through the lower bound on the scalar curvature. The main characteristic of these relations is that they show how the local properties in some regions affect the local properties in regions far away.   
 
The main inequality, which is introduced in Section \ref{ss1}, is a slight generalization of a result that to our knowledge first appeared in the work of Castillon \cite{Castillon}. The formula provides suitable geometric information when the second variation operator is used with radial functions. Its proof is postponed until the last section. In Section \ref{SC} we present the first applications of the main inequality. The main goal is to obtain suitable differential inequalities for the areas of metric balls. After integration, the equations establish a non-local relation between the areas of balls of different radius. As is explained in Section \ref{SC} these applications have interest in themselves but in Section \ref{S3} we utilize them to get more elaborated results. 
In the first of them we show that when the stable minimal surface is the two sphere and the ambient scalar curvature is non-negative then the first non-zero eigenvalue of the Laplacian (on the surface) is controlled from above and below by one over the diameter square of the surface. This complements the Hersh inequality \cite{Hersch} for the first non-zero eigenvalue of the Laplacian of two-spheres. In the second application we show that the well known isoperimetric properties that hold on Riemannian surfaces with Gaussian curvature bounded below also hold on stable minimal surfaces lying inside three-manifolds with scalar curvature bounded below. In particular we show non-collapse at a finite distance from a non-collapsed region (see later for a precise definition). This allows us to obtain a compactness results for families of pointed stable minimal surfaces with $L^{p}$ ($p>1$) curvature uniformly bounded above. 
In Section \ref{S3.3} we prove positively a question raised by Fischer-Colbrie and Schoen in \cite{Fischer-Colbrie-Schoen}. Precisely we prove that if the ambient scalar curvature is non-negative then any complete and stable minimal surface diffeomorphic to a punctured disc (``a tube" for us) is flat and totally geodesic. Thus, according to the Fischer-Colbrie/Schoen classification \cite{Fischer-Colbrie-Schoen} if a complete non-compact stable minimal surface is non-flat then it has to be topologically the plane. In Section \ref{S3.4} we illustrate how to use the geometric relations expressed in Theorem \ref{T1} to understand the notion of local isoperimetric collapse on three-manifolds with non-negative scalar curvature. This notion is crucial to make further progress  in the theory of convergence/collapse of Riemannian three-manifolds having a priori bounds on the $L^{p}$ ($p>3/2$) curvature (for a treatment of this subject and potential applications to General Relativity see \cite{Reiris}). 
 
\vs
\vs
{\center \subsection{Notation.}}
\vs
$(\Su,h)$ will be a complete orientable two-surface (compact or not and possible with smooth boundary) immersed inside an orientable  Riemannian 
three-dimensional manifold $(\Sigma,g)$, where $h$ is the induced metric (from now on we will omit writing down $h$ and $g$ explicitly). The scalar curvature of the metric $g$ will be denoted by $R$ and the Gaussian curvature of the metric $h$ by $\kappa$.

Let $U$ be a two-manifold, possibly with smooth boundary. We define the {\it Radius} of $U$ as $Rad(U)=\sup\{dist(p,\partial U),\ p\in U\}$. The {\it Radius centered at a set $C$} (for instance a point $C=\{q\}$) is defined by $Rad(C,U)=\sup\{dist(p,C),\ p\in U\}$. In this notation $Rad(U)=Rad(\partial U,U)$.  Also  {\it balls
of center a set $C$ and radius $r$} are denoted as $B(C,r)$ ($B(C,r)=\{p\in S/\ dist(p,C)\leq r\}$). 

{\center \section{The main inequality.\label{S2}}}
{\center \subsection{Using the stability inequality with radial functions.\label{ss1}}}

We will assume throughout that the scalar curvature $R$ of $g$ is uniformly bounded below by $R_{0}$. In practice one can take $R_{0}=\inf\{R(p),p\in \Su\}$ in each of the statements that follow. We will assume too that $\Su$ is stable, namely that for every $f$ of compact support in the interior of $\Su$, we have the inequality
\ben
\int_{\Su}|\nabla f|^{2}+\kappa f^{2}dA\geq \int_{\Su}(\frac{R}{2}+\frac{|H|^{2}}{2})f^{2}dA,
\een

\n where $H$ is the second fundamental form of $\Su$ in $\Sigma$ and $\kappa$ the Gaussian curvature.  We then have the inequality 
\be\label{svf}
\int_{\Su}|\nabla f|^{2}+\kappa f^{2}dA\geq \int_{\Su}\frac{R_{0}}{2}f^{2}dA.
\ee

The next statement explains a geometric inequality (deduced from (\ref{svf})) that gives us information on the geometry of such surfaces. The remarkable fact is that the inequality depends on the ambient curvature only through the lower bound on the scalar curvature.  This fact makes the technique uniquely suited to study the geometry of three-manifolds (with scalar curvature uniformly bounded below) through stable minimal surfaces.

\vs
\begin{Theorem}\label{T1} Let $\Su$ be a compact surface with smooth boundary consisting of a finite set of closed curves $\partial \Su=\{\ell_{1}\cup\ldots\cup \ell_{n}\}$. Consider the distance function $r(p)=dist(p,\ell_{i_{1}}\cup\ldots\cup \ell_{i_{j}})$ where $j $ is a fixed number between $1$ and $n$. Fix a distance $L$ less or equal than the distance between $\ell_{i_{1}}\cup\ldots\cup \ell_{i_{j}}$ and $\partial \Su\setminus (\ell_{i_{1}}\cup\ldots\cup \ell_{i_{j}})$ and consider the radial function $f(p)=1-\frac{r(p)}{L}$ if $r(p)\leq L$ and zero otherwise. Then we have
\be\label{svc}
\int_{\Su}|\nabla f|^{2}+\kappa f^{2}dA\leq 2\frac{l}{L}+l'-\frac{A}{L^{2}},
\ee

\n were $l=length(\ell_{i_{1}})+\ldots +length(\ell_{i_{j}})$ and $l'=length(\ell_{i_{1}})'+\ldots+length(\ell_{i_{j}})'$ where $'$ denotes the derivative when the curve is displaced in its inward normal direction. Also $A=Area(B(\ell_{i_{1}}\cup\ldots\cup \ell_{i_{j}},L))$.
\end{Theorem}

\vspace{0.2cm}
Theorem \ref{T1} is a mild generalization of a result that to our knowledge first appeared in the work of Castillon \cite{Castillon} (Castillon considers the distance function to a point). Similar results to Castillon's where subsequently used by Rosenberg and Espinar \cite{Rosenberg-Espinar} in topics somehow related with the present article. Their proof (as provided in \cite{Castillon} and \cite{Rosenberg-Espinar}) follows after an integration by parts in the second variation formula (Pogorelov \cite{Pogorelov}). Such integration by parts was used by Colding and Minicozzi \cite{Colding-Minicozzi} for stable minimal surfaces inside manifolds having scalar curvature uniformly bounded below (the trick however was used by other authors as well). The crucial difference made in \cite{Castillon} lies in the fact that the integration by parts is carried out beyond the locus (see later). To achieve this, a non-obvious formula due to Shiohama and Tanaka \cite{Shiohama-Tanaka} is used. For the sake of completeness we will provide a proof of Theorem \ref{T1} at the end of the article. It will be carried out assuming a minor technical simplification that we believe makes the proof more conceptual and explicit. Every step  can be tracked down to Pogorelov's integration and the results of Shiohama and Tanaka in \cite{Shiohama-Tanaka}.         

\vs
{\center \subsection{Some consequences of Theorem \ref{T1}.}\label{SC}}
\vs

We present now a set of propositions that will be useful in the applications. We remark however that Theorem \ref{T1} can applied in a great number of situations and in different ways. 
From now on, and to simplify the analysis, assume that $\Su$ is a stable minimal surface that is either compact, or non-compact but without boundary. We will use the terminology {\it numeric} to mean a number independent on any aspect of the hypothesis. 

\vspace{0.2cm}
\begin{Proposition}\label{P1} Assume $R\geq 0$. There is a numeric $c>0$ such that for any $L'<L\leq Rad({p},\Su)$ it is
\be
A(B(p,L'))\geq c(\frac{L'}{L})^{2}A(B(p,L)).
\ee

\end{Proposition}

This statement is a first indication of the important fact that for stable minimal surfaces lying inside three-manifolds with scalar curvature uniformly bounded below, the local geometry is strongly entangled with the global geometry. In this case what Proposition \ref{P1} says is that the area of balls is definitely influenced by the area of bigger balls (with the same center). 

\vspace{0.2cm}
\n {\bf Proof:} 

For any $r<L'/2$ consider the set ${\mathcal{B}}=\Su\setminus \overline{B(p,r)}$. If a connected component (denoted by c.c. below) $U$ of ${\mathcal{B}}$ has $Rad(U)\leq L'-r$ then we join that component to $\overline{B(p,r)}$, the resulting set will be denoted by ${\mathcal{C}}$. (We will assume below that $\partial {\mathcal{C}}$ is smooth, if not, one can get the same formulae below by smoothing and taking a limit back). 

We chose now the function $f$ to be used in the operator
\be\label{svo}
f\rightarrow \int_{\Su}|\nabla f|^{2}+\kappa f^{2}dA.
\ee

\n Define

\begin{equation*}
\begin{array}{ll}
f(p)=1-\frac{dist(p,\partial {\mathcal{C}})}{Rad({\mathcal{C}})} & {\rm if}\ p\in {\mathcal{C}},\\
&\\
f(p)=1-\frac{d(p,\partial U)}{Rad(U)} & {\rm if}\ p\in U,\ {\rm and\ U\ a\ c.c.\ of\ } {\mathcal{B}}\ {\rm with}\ L-r\geq Rad(U)>L'-r,\\
&\\  
f(p)=1-\frac{dist(p,\partial U)}{L-r} & \ {\rm if}\ p\in U,\ {\rm U\ a\ c.c.\ of}\ {\mathcal{B}}\ {\rm with}\ Rad(U)>L-r,\ {\rm and}\ dist(p,\partial U)\leq L-r,\\
&\\
f(p)=0 & {\rm otherwise.}
\end{array}
\end{equation*}

\n Denote $l_{k}=length(\partial U_{k})$ for every component $U_{k}$, $k=1,\ldots, m$ of ${\mathcal{B}}$ with $L-r\geq Rad(U_{k})> L'-r$. Also denote $\bar{l}_{k}=length(\partial U_{k})$ where $U_{k}$, $k=m+1,\ldots,\bar{m}$ are the connected components of ${\mathcal{B}}$ with $Rad(U_{k})>L-r$. Finally denote $l=\sum_{k=1}^{k=m}l_{k}+\sum_{k=m+1}^{k=\bar{m}}\bar{l}_{k}$ the length of the full boundary of ${\mathcal{C}}$. Theorem \ref{T1} gives (note the cancelation of the $l'$-type of terms)
\be\label{iP1}
2\frac{l}{Rad({\mathcal{C}})}+2\sum_{k=1}^{k=m}\frac{l_{k}}{Rad(U_{k})}+2\sum_{k=m+1}^{k=\bar{m}}\frac{\bar{l}_{k}}{L-r}\geq \frac{A({\mathcal{C}})}{Rad({\mathcal{C}})^{2}}+\sum_{k=1}^{k=m}\frac{A(U_{k})}{Rad(U_{k})^{2}}+\sum_{k=m+1}^{k=\bar{m}} \frac{A(U_{k})}{(L-r)^{2}}.
\ee

\n This formula is useful in itself, however in order to apply it to prove the statement, we need to simplify it further. We treat the left hand side first. Observe that $Rad({\mathcal{C}})$, $Rad(U_{k})$ for $k=1,\ldots,m$ and $L-r$ are greater or equal than $r$. Therefore the left hand side is less or equal than $4l/r$. To simplify the right hand side note that any two $q_{1}$ and $q_{2}$ in $\partial B(p,r)$ are connected through two radial geodesics meeting at $p$ and therefore we have $Rad({\mathcal{C}})\leq L'-r+2r=L'+r\leq (3/2)L'\leq (3/2)L$. Note also that $Rad(U_{k})\leq L-r\leq L$ for $k=1,\ldots,m$ and that $L-r\leq L$. Therefore we get that the right hand side is greater or equal than $(2/(3L))^{2}A(B(p,L))$. Finally note that $A(B(p,r))'\geq l$. Thus we get
\ben
\frac{A'}{r}\geq \frac{1}{9}\frac{A(B(p,L))}{L^{2}}.
\een

\n Integrating from $0$ to $r$ gives 
\ben
A(r)\geq \frac{1}{18}\frac{A(L)}{L^{2}}r^{2}.
\een

\n Evaluating at $r=L'/2$ and noting that $A(L')\geq A(L'/2)$ we get
\ben
A(L')\geq \frac{1}{72}\frac{L'^{2}}{L^{2}}A(L).
\een

\ep

\vspace{0.2cm}
The well known inequality $A(B(p,L))\leq 2\pi L^{2}$ (use $f(o)=1-dist(o,p)/L$ if $dist(o,p)\leq L$ and zero otherwise in the operator (\ref{svo})) can be complemented a bit in some situations using the previous proposition. One example is the following.

\begin{Corollary}\label{C1} Assume $R\geq 0$. Let $\Su$ be a complete non-compact stable minimal surface. Then if at a point $p$ it is
\ben
\limsup_{L\rightarrow \infty}\frac{A(B(p,L))}{L^{2}}>0,
\een

\n then 
\ben
\liminf_{L\rightarrow \infty}\frac{A(B(p,L))}{L^{2}}>0.
\een

\end{Corollary}

\vspace{0.2cm}
Proposition \ref{P1} does not tell how the the area of metric balls grow if we do not know a priori the area of a certain ball. The following Proposition remedies in part this problem.

\begin{Proposition}\label{P3} Assume $R\geq 0$. Let $p$ and $o$ be two points in $\Su$. Make $A(r)=Area(B(p,r))$. Then if $0<r< dist(p,o)$ we have 
\be\label{iP3}
2 A'(\frac{1}{r}+\frac{1}{dist(p,o)-r})\geq \frac{A}{4r^{2}}+\frac{A_{*}}{(dist(p,o)-r)^{2}},
\ee

\n where $A_{*}$ is the area of the intersection of the connected component of $\Su-\overline{B(p,r)}$ containing $o$ with $B(p,dist(p,o))$. 
\end{Proposition}

\n {\bf Proof:} 

Let ${\mathcal{L}}$ be the boundary of the connected component of $\Su-\overline{B(p,r)}$ containing $o$. To simplify the notation below denote such component $U$, also denote 
$D=Rad({\mathcal{L}},B(p,r))$. In the operator (\ref{svo}) consider the function $f$ 
\ben
\begin{array}{ll}
f(q)=1-\frac{dist(q,{\mathcal{L}})}{D} &\ {\rm if\ } q\in \Su\setminus U,\ {\rm and}\ dist(q,{\mathcal{L}})\leq D,\\
&\\
f(q)=1-\frac{dist(q,{\mathcal{L}})}{dist(p,o)-r}&\ {\rm if}\ q\in U,\ {\rm and}\ dist(q,{\mathcal{L}})\leq dist(p,o)-r,\\
&\\
f(q)=0&\ {\rm otherwise.}
\end{array}
\een

\n From Theorem \ref{T1} we get
\ben
2length({\mathcal{L}})(\frac{1}{D}+\frac{1}{dist(p,o)-r})\geq \frac{A^{*}}{D^{2}}+\frac{A_{*}}{(dist(p,o)-r)^{2}},
\een

\n where $A^{*}$ is the area of the set of points in $\Su\setminus U$ at a distance at most $D$ from ${\mathcal{L}}$. As we have, $r\leq D\leq 2r$, $A^{*}\geq A$ and $A'\geq length({\mathcal{L}})$ the result follows.\ep

\vspace{0.2cm}
We get the following immediate Corollary complementing Corollary \ref{C1}.
\vs
\begin{Corollary} Assume $R\geq 0$. Let $\Su$ be a non-compact stable minimal surface. For any $p$ in $\Su$ there is $c(p)$ such that
\ben
A(B(p,r))\geq c(p)r^{\frac{1}{8}}.
\een 

\end{Corollary}
\vs
\n {\bf Proof:} 

Picking $o$ at infinity in Proposition \ref{P3} we get the inequality (make $A(r)=A(B(p,r))$
\ben
\frac{A'}{A}\geq \frac{1}{8}\frac{1}{r}.
\een

\n Integrating between $r_{0}=1$ and $r>1$ we get the result.\ep
\vs

Another application of inequality (\ref{iP3}) is an isoperimetric inequality for metric balls, precisely we have

\begin{Corollary} Assume $R\geq 0$. Let $p$ and $o$ be points in $\Su$. Make $A(r)=Area(B(p,r))$ and $l=length(\partial B(p,r))$. There is a numeric $c>0$ such that if $0<r<dist(p,o)/2$  
\ben
\frac{l(r)}{A^{\frac{1}{2}}(r)}\geq c\frac{A(B(o,dist(p,o)/2))^{\frac{1}{2}}}{dist(p,o)}.
\een
\end{Corollary}

\n {\bf Proof:} 

From inequality (\ref{P1}) we have
\ben
2\frac{l(r)}{r}\geq \frac{A(r)}{4r^{2}},
\een

\n and therefore
\be\label{c1}
\frac{l}{A^{\frac{1}{2}}}\geq \frac{1}{8}\frac{A^{\frac{1}{2}}}{r}.
\ee

\n We also have from Proposition \ref{P1} (take $L=(3/2)dist(p,o)$ and $L'=r$)
\ben
A(r)\geq c\frac{r^{2}}{dist(p,o)^{2}}A(B(o,dist(p,o)/2)),
\een

\n and thus
\be\label{c2}
\frac{A^{\frac{1}{2}}}{r}\geq (c\frac{A(B(o,dist(p,o)/2))}{dist(p,o)^{2}})^{\frac{1}{2}}.
\ee

\n Combining (\ref{c1}) and (\ref{c2}) the result follows.\ep 

\vs
A similar result than Proposition \ref{P1} holds when $R_{0}<0$ however it does not for balls $B(p,L')$ and $B(p,L)$ of arbitrary radiuses. If $R_{0}<0$ we want to exploit the non-negativity of the operator
\ben
f\rightarrow \int_{\Su}|\nabla f|^{2}+\kappa f^{2} +\frac{R_{0}}{2}f^{2} dA.
\een

\begin{Proposition}\label{P33} Suppose $R_{0}<0$. Let $L'<L\leq Rad(p,\Su)$ with $-9 R_{0}L^{2}\leq 1$. Then there is a numeric $c$ such that 
\ben
A(B(p,L'))\geq c(\frac{L'}{L})^{2} A(B(p,L)).
\een

\end{Proposition}
\n {\bf Proof:}
The proof is a modification of the argument of Proposition \ref{P1}. Indeed, following the same notation, we get a similar inequality than (\ref{iP1}) but this time we gain the extra term
\ben
\frac{R_{0}}{2}(A({\mathcal{C}})+\sum_{k=1}^{k=\bar{m}}A(U_{k})),
\een

\n on its right hand side. Thus, using the condition $-9 R_{0}L^{2}\leq 1$ we get
\ben
4\frac{l}{r}\geq ((\frac{2}{3})^{2}\frac{1}{L^{2}}+\frac{R_{0}}{2})A({\mathcal{C}})+\sum_{k=1}^{k=\bar{m}}(\frac{1}{L^{2}}+\frac{R_{0}}{2})A(U_{k})\geq \frac{1}{3L^{2}}A(B(p,L)),
\een

\n and the result follows.\ep

\vs
{\center \section{Applications.\label{S3}}}
{\center \subsection{Controlling the first non-zero eigenvalue of the Laplacian from the diameter.\label{S3.1}}}
\vs

\begin{Proposition}\label{P5} Suppose $R\geq 0$ and that $\Su$ is a stable minimal surface diffeomorphic to the two-sphere. Let $\lambda>0$ be the first non-zero eigenvalue of minus the standard Laplacian. Then there are numeric $c_{1}>0$ and $c_{2}>0$ such that
\be\label{iP5}
c_{1}\frac{1}{diam(\Su)^{2}}\leq \lambda \leq c_{2}\frac{1}{diam(\Su)^{2}}.
\ee
\end{Proposition}

\n {\bf Proof:}

We estimate first the Cheeger constant \cite{Chavel} (pg. 109) that we will 
denote by $\xi$. Recall that 
\ben
\xi=\inf\frac{length(\ell)}{\min\{A_{1},A_{2}\}},
\een

\n where $\ell$ ranges over all closed loops in $\Su$ and $A_{1}$ and $A_{2}$ are the areas of the connected components $U_{1}$ and $U_{2}$ of
$\Su\setminus \ell$. Denote $L_{1}=Rad(\ell,U_{1})$ and $L_{2}=Rad(\ell,U_{2})$. Applying Theorem \ref{T1} one has (for a certain sequence of loops $\{\ell_{i}\}$)
\ben
\xi=\lim \frac{\ell_{i}}{\min\{A_{1,i},A_{2,i}\}}\geq \frac{1}{2}\frac{1}{\min\{A_{1,i},A_{2,i}\}}\frac{A_{1,i}/L_{1,i}^{2}+A_{2,i}/L_{2,i}^{2}}{1/L_{1,i}+1/L_{2,i}}\geq \frac{1}{2}\frac{1/L_{1,i}^{2}+1/L_{2,i}^{2}}{1/L_{1,i}+1/L_{2,i}},
\een
\ben
\geq \frac{1}{2}\min\{1/L_{1,i},1/L_{2,i}\}\geq \frac{1}{2}\frac{1}{diam(\Su)}.
\een

\n It is well known on the other hand \cite{Chavel} that $\lambda\geq \frac{\xi}{4}$, thus showing one direction of the inequality 
(\ref{iP5}). We prove now the other direction of the inequality. Pick $s$ and $o$ such that $dist(s,o)=diam(\Su)$. Consider the function $f(p)=1-\frac{dist(p,o)}{diam(\Su)}$. It is
\ben
\lambda\leq \frac{\int_{\Su}|\nabla f|^{2}d A}{\int_{\Su} |f-\bar{f}|^{2} dA},
\een

\n where $\bar{f}=(\int_{\Su}f dA)/A$ is the average of $f$. By Proposition \ref{P1} the balls $B(s,diam/3)$ and $B(o,diam/3)$ carry at least a definite fraction of the total area $A$. It follows that $\bar{f}\leq c<1$ where $c$ is numeric. We have
\be\label{i2P5}
\int_{\Su}|f-\bar{f}|^{2}dA\geq \int_{\{1-r/d\geq (1+c)/2\}}(\frac{1-c}{2})^{2}dA\geq (1/c_{2})A(\Su),
\ee

\n where the last inequality follows from Proposition \ref{P1} and $c_{2}$ is numeric. The inequality (\ref{i2P5})
now gives
\ben
\lambda\leq c_{2}\frac{A(\Su)}{diam(\Su)^{2}A(\Su)}=c_{2}\frac{1}{diam(\Su)^{2}},
\een

\n finishing the proof of the inequality (\ref{iP5}).\ep 

\vs
{\center \subsection{Non-collapse at a fine distance from a non-collapsed region.}}
\vs

Let $\Su$ be an arbitrary stable surface that to simplify the analysis we chose either compact or non-compact but without boundary. We start this section by recalling the well known definition of volume radius. Given $\delta>1$ define 
\ben
\vru(p)^{\delta}=\sup\{r/Rad(p,\Su)\geq r>0\ and\ A(B(p,r))\leq \pi\delta^{2}r^{2}\},
\een
\ben
\vrd(p)^{\delta}=\sup\{r/Rad(p,\Su)\geq r>0\ and\ A(B(p,r))\geq \frac{\pi}{\delta^{2}}r^{2}\}.
\een

\n Observe that if $\delta_{1}\geq \max\{\delta_{2},\delta_{3}\}$ then 
\ben
\vru(p)^{\delta_{1}}\geq \max\{\vru(p)^{\delta_{2}},\vru(p)^{\delta_{3}}\},
\een
\ben
\vrd(p)^{\delta_{1}}\geq \max\{\vrd(p)^{\delta_{2}},\vrd(p)^{\delta_{3}}\}.
\een

\n In particular if $\vrd(p)^{\delta_{1}}\geq a_{1}$ and $\vrd^{\delta_{2}}(p)\geq a_{2}$ then $\vrd(p)^{\delta_{1}+\delta_{2}}\geq \max\{a_{1},a_{2}\}$ and similarly for $\vru^{\delta}$. These observations will be used later.

Define the volume radius $\nu(p)^{\delta}$ at $p$ as
\ben
\nu(p)^{\delta}=\sup\{r/Rad(p,\Su)\geq r>0\ and\ \frac{\pi}{\delta^{2}}r^{2}\leq A(B(p',r))\leq \pi\delta^{2}r^{2},\forall B(p',r)\subset B(p,\nu(p)^{\delta})\}.
\een

\n It is clear that estimates on $\vru(p)^{\delta}$ and $\vrd(p)^{\delta}$ around a point $o$ would give estimates on $\nu(o)^{\delta}$. The quantities $\vru$ and $\vrd$ are easier to manipulate and for this reason we use them in the next theorem (instead of $\nu$).  

The volume radius is a basic local invariant. As is well known \cite{Anderson} an upper bound on $1/\nu(o)^{\delta}$ and $\|\kappa\|_{L^{p}(B(o,\nu(o)^{\delta}))}$ ($p>1$) gives control on the $H^{2,p}_{\{x_{i}\}}$-Sobolev norm of the metric entrances $h_{ij}-\delta_{ij}$ on a harmonic coordinate system $\{x_{i}\}$ covering a ball $B(o,r)$ where $r$ is also controlled from below. The analysis that we will carry below on the volume radius combined with this standard fact will give the compactness result of Proposition \ref{C}.

We pass now to study the volume radius on stable surfaces immersed on three-manifolds with $R\geq R_{0}$. Let us observe first that for any point $p$ in $\Su$ and for any $r\leq Rad(p,\Su)$ we have
\ben
A(B(p,r))\leq \frac{2\pi r^{2}}{(1+\frac{\tilde{R_{0}}r^{2}}{2})}.
\een

\n where $\tilde{R}_{0}=inf\{R_{0},0\}$ and also $r^{2}R_{0}/2\geq -1$ (this known property follows after taking the trial function, $f(q)=1-dist(q,p)/r$ if $dist(q,p)\leq r$ and zero otherwise, in the stability inequality). It follows that if $R_{0}\geq 0$ then 
\ben
\vru(p)^{2}\geq \sup\{d(p,o),o\in \Su\}\geq \frac{diam(\Su)}{2},
\een

\n (independently on whether $diam=\infty$ or not). Similarly if $R_{0}<0$ we have 
\ben
\vru(p)^{2}\geq \min\{\frac{diam(\Su)}{2},\frac{1}{\sqrt{-R_{0}}}\}.
\een

\n For this reason we will occupy ourselves in the next Theorem with an estimate for $\vrd(p)^{\delta}$, and for arbitrary points $p$ in $\Su$. We will denote numeric constants below by $c_{i},\ i=1,2,\ldots$. If used in different items we do not mean that necessarily they are the same.

\vs
\begin{Theorem}\label{T2} Let $\Su$ be a complete and stable minimal surface that (to simplify the treatment) we assume without boundary. Let $p$ and $o$ be arbitrary points in $\Su$. Make $\nu(o)^{\delta}=\nu_{0}$.
\begin{enumerate}
\item If $R_{0}>0$ then 
\be\label{vr(i)}
\underline{\nu}(p)^{c_{1}\delta(\frac{dist(p,o)}{\nu_{0}}+1)}\geq \max\{\frac{\nu_{0}}{2},dist(p,o)\}.
\ee

\n Moreover $\Su$ is compact and we have the estimate $diam(\Su)\leq c_{2}\frac{\delta}{R_{0}\nu_{0}}$.

\item If $R_{0}=0$, then we have the same estimate as before for $\underline{\nu}^{\bar{\delta}}(p)$ (with the same $\bar{\delta}$). In addition if the total area $A(\Su)$ is finite then
$\Su$ is compact and we have the estimate 
\be\label{*}
diam(\Su)\leq c_{1}\frac{A(\Su)^{8}\delta^{16}\nu_{0}}{\nu_{0}^{16}}.
\ee

\n If instead the total area is infinite, then the area of $B(o,r)$ grows like 
\be\label{**}
A(r)\geq c_{2} (\frac{\nu_{0}}{\delta})^{2}(\frac{r}{\nu_{0}})^{\frac{1}{8}}.
\ee

\item If $R_{0}<0$, then
\be\label{vr(iii)}
\underline{\nu}(p)^{c_{1}\delta 2^{(c_{2}\frac{dist(p,o)}{\min\{\frac{c_{3}}{\sqrt{-R_{0}}},c_{4}\nu_{0}\}} )}}\geq \min\{\frac{1}{9\sqrt{-R_{0}}},\frac{\nu_{0}}{9}\}.
\ee
\end{enumerate}
\ep\end{Theorem}

A comment on this statement is in order. We say that a family of surfaces has the {\it non-collapse at a finite distance} property iff {\it there are $\delta_{1}(\nu_{0},D)>1$ and $\nu_{1}(\nu_{0},D)>0$ such that if $\nu(o)^{\delta}\geq \nu_{0}$ and $dist(p,o)\leq D$ then $\nu(p)^{\delta_{1}}\geq \nu_{1}$}. What the Theorem above shows is that the family of stable minimal surfaces immersed on three-manifolds with scalar curvature uniformly bounded below ($R\geq R_{0}$ ($R_{0}$ fixed)) enjoys the {\it non-collapse at a finite distance property}.   

\vs
\n {\bf Proof:}

{\it Item 1}. Suppose $d\geq \frac{\nu_{0}}{2}$ then by Proposition \ref{P1} we have (make $A(r)=A(B(p,r))$ and $d=dist(p,o)$)
\ben
A(r)\geq c_{1}(\frac{r}{d})^{2}A(B(p,d))\geq c_{2}(\frac{r\nu_{0}}{d\delta})^{2},
\een

\n as long as $r\leq d$. Thus if $d\geq \nu_{0}/2$ then there is $c_{1}$ such that
\ben
\vrd^{c_{1}\delta \frac{d}{\nu_{0}}}(p)\geq dist(p,o)=\max\{\frac{\nu_{0}}{2},dist(p,o)\}.
\een

\n On the other hand if $d\leq \nu_{0}/2$ we have $A(r)\geq \pi r^{2}/\delta^{2}$ whenever $r\leq \nu_{0}/2$. Thus if $d\leq \nu_{0}/2$ then there is $c_{2}$ such that 
\ben
\vru(p)^{c_{2}\delta}\geq \frac{\nu_{0}}{2}=\max\{\frac{\nu_{0}}{2},dist(p,o)\}.
\een

\n Combining this with the observation at the beginning of the section we deduce that there is $c_{3}$ such that 
\ben
\underline{\nu}(p)^{c_{3}\delta(\frac{d}{\nu_{0}}+1)}\geq \max\{\frac{\nu_{0}}{2},dist(p,o)\},
\een

\n as claimed. Let us prove now the second part of the {\it Item 1}. Let $\gamma$ be a length minimizing geodesic joining $p$ with $o$. Let $q$ be the middle point between $p$ and $o$. By Proposition \ref{P1} we have $A(B(q,d/6))\geq c_{1} A(B(q,d/2))\geq c_{2}(\frac{\nu_{0}}{\delta})^{2}$. Consider now the stability inequality (\ref{svf}). We will choose the function $f$ as we did in Proposition \ref{P1} but choosing $L=d$, $L'=d/3$ and $r< d/6$. Following the same terminology we note that $B(q,d/6)$ belongs to a connected component $U$ with $Rad(U)\geq L-r\geq 5d/6$. We note next that if $t\in B(q,d/6)$ then $dist(t,\partial {\mathcal{C}})\leq 2d/6+d/3=2d/3$ and therefore there is a numeric $c_{3}>0$ which is a lower bound for $f^{2}/2$ over $B(q,d/6)$ for any $r<d/6$. Thus there is a lower bound for the right hand side of (\ref{svf}) of the form $c_{3}R_{0}A(B(q,d/6))$, which is greater or equal than $c_{4}R_{0}(\nu_{0}/\delta)^{2}$, while for the left hand side of (\ref{svf}) there is the upper bound $4A'/r$ ($r\leq d/6$). Therefore we get 
\ben
A(r)\geq c_{5} r^{2}(\nu_{0}/\delta)^{2}R_{0}.
\een

\n On the other hand (taking $f(p')=1-\frac{dist(p',p)}{dist(p,o)}$) we have 
\ben
2\pi\geq c_{6}R_{0}A(B(p,d/6)).
\een
     
\n The two inequalities above imply (take $r=d/6$) that for any two points $p$ and $o$ we have
\ben
dist(p,o)\leq c_{7}\frac{\delta}{R_{0}\nu_{0}},
\een

\n and the result follows. 

{\it Item 2}. We need to prove only the second part of the statement {\it item 2}. Recall that from Proposition \ref{P3} we have (make $A(r)=A(B(o,r))$)
\ben
A'(\frac{Rad(o)}{r(Rad(o)-r)})\geq \frac{1}{8}\frac{A}{r^{2}}.
\een

\n Thus
\ben
\frac{A'}{A}\geq \frac{1}{8}(\frac{1}{r}-\frac{1}{Rad(o)}).
\een

\n Integrating we get
\ben
A(r)\geq A(\nu_{0})(\frac{r}{\nu_{0}})^{\frac{1}{8}}e^{-\frac{r-\nu_{0}}{8Rad(o)}}.
\een

\n If $A(\Su)$ is finite then by taking $r=Rad(o)\geq diam(\Su)/2$ we get inequality (\ref{*}). On the other hand if the total area of $\Su$ is infinite then necessarily $Rad(o)=\infty$ and the equation (\ref{**}) follows.

{\it Item 3}. Let $\gamma$ be a length minimizing geodesic joining $p$ and $o$. Let $d_{0}=\min\{\frac{1}{3\sqrt{-R_{0}}},\nu_{0}/3\}$. Suppose first that $dist(p,o)\geq d_{0}$. Consider now the following set of disjoint balls: the ball of center at a distance of $d_{0}/3$ along $\gamma$ from $o$ and radius $d_{0}/3$; the ball of center at a distance $d_{0}$ from $o$ and radius $d_{0}/3$; the ball of center at a distance $5d_{0}/3$ from $o$ and radius $d_{0}/3$ and so on. We have $n=[3d(p,o)/(2d_{0})]$ of such balls between $p$ and $o$.  According to Proposition \ref{P33}
we have
\ben
A_{k}=\frac{1}{c}A_{k-1},
\een

\n where $A_{k}$ is the ares of the $k-$ball. Thus
\ben
A_{k}\geq \frac{A_{1}}{c^{k-1}},
\een

\n for $k=2,3,\ldots$ until $k=n$.  Therefore we have
\ben
A_{k}\geq \frac{\pi d_{0}^{2}}{9c^{k-1}\delta^{2}}.
\een

\n Note that because $dist(p,o)\geq d_{0}\leq \frac{\nu_{0}}{9}$ we have $n\geq 2$. Now $p$ is at a distance between $d_{0}$ and $5d_{0}/3$ from the center
of the $(n-1)$-ball. Again Proposition \ref{P33} gives 
\be\label{AC}
A(B(p,r))\geq c_{1}\frac{A_{n-1}}{d_{0}^{2}}r^{2}\geq \frac{c_{1}\pi}{9c^{n-2}\delta^{2}}r^{2}\geq \frac{c_{2}\pi}{\delta^{2}2^{c_{3}\frac{dist(p,o)}{d_{0}}}}r^{2},
\ee

\n valid for $r\leq d_{0}/3$. Noting that if in (\ref{AC}) we have $c_{2}\geq 1$ then one can take $c_{2}=1$ naturally preserving the inequality and thus having 
\ben
A(B(p,r))\geq \frac{c_{2}\pi r^{2}}{\delta^{2}2^{c_{3}\frac{d}{d_{0}}}},
\een

\n with $(2^{c_{3}d/\nu_{0}})/c_{2}\geq 1$. In case $dist(p,o)\leq d_{0}\leq \frac{\nu_{0}}{3}$ we have $A(B(p,r))\geq \frac{\pi}{\delta^{2}}r^{2}$, and valid
for $r\leq \frac{\nu_{0}}{3}$ and therefore for $r\leq d_{0}/3$. The inequality (\ref{vr(iii)}) now follows directly from the observation at the beginning of the section.\ep

\vs
It is direct to show, using the standard volume (area) comparison, that the volume radius in Riemannian surfaces having Gaussian curvature $\kappa$ bounded below by $\kappa_{0}$ has the following properties.
(We use below the terminology: {\it a positive quantitie $A$ is controled by a positive quantity $B$ iff for any $B_{0}$ there is $A_{0}$ such that if $B\leq B_{0}$ then $A\leq A_{0}$. The same definition applies if $A$ and $B$ are a set of positive quantities.})  

\begin{Theorem} Let $S$ be a complete Riemannian surface with Gaussian curvature $\kappa\geq\kappa_{0}$.

\begin{enumerate}
\item If $\kappa_{0}>0$ then $S$ is compact and $diam(S)$ is controlled by $1/\kappa_{0}$. Moreover $\Gamma$ and $1/\nu(p)^{\Gamma}$ are controlled by $1/\nu(o)^{\delta}$.

\item If $\kappa_{0}=0$ then $\Gamma$ and $1/\nu(p)^{\Gamma}$ are controlled by $1/\nu(o)^{\delta}$ and $dist(p,o)$. Moreover

\begin{enumerate}
\item If $S$ is compact then $diam(S)$ is controlled by $Area(S)$ and $1/\nu(o)^{\delta}$.
\item If $S$ is non-compact then $Area(S)=\infty$.
\end{enumerate}

\item If $\kappa_{0}<0$ then $\Gamma$ and $\nu(p)^{\Gamma}$ are controlled by $|\kappa_{0}|$, $1/\nu(o)^{\delta}$ and $dist(p,o)$. 
\end{enumerate}
\end{Theorem}  
  
\n A direct inspection shows that these are the same properties that have shown to hold on stable minimal surfaces if we make the correspondence $\kappa_{0} \longleftrightarrow R_{0}$. From some point of view at least, the family of stable minimal surfaces on three-manifolds with $R\geq R_{0}$ and the family of Riemannian surfaces with Gaussian curvature $\kappa\geq \kappa_{0}$ display and interesting parallel.

A distinguished Corollary of Theorem \ref{T2} is the following compactness result. We state it in a simple form. For a proof see \cite{Anderson}.
\begin{Proposition}\label{C}
Consider the family ${\mathcal{F}}_{R_{0},\Lambda_{0},\nu_{0},p}$ ($p>1$) of complete pointed $H^{3,p}$-Riemannian surfaces $(\Su,h,o)$ without boundary, satisfying:
\begin{enumerate}
\item The stability inequality (\ref{svf}),
\item $\|\kappa\|_{L^{p}_{h}(\Su)}\leq \Lambda_{0}$,
\item $\nu(o)^{\delta}\geq \nu_{0}$.
\end{enumerate}

Then the family is sequentially compact in the weak $H^{2,p}$-topology.
\end{Proposition}    

\n Above $H^{k,p}$ represents (generically) the Sobolev space of distributions with $k$-derivatives in $L^{p}$. Thus fixed a chart $\{x\}$, the transition functions to other charts in $\Su$ are in $H^{3,p}_{\{x\}}$ and $h_{ij}$ is in $H^{2,p}_{\{x\}}$.

{\center \subsection{If $R_{0}\geq 0$ then complete and stable tubes are flat (answering a question of Fischer-Colbrie and Schoen).}\label{S3.3}}
\vs
\begin{Lemma}\label{L2} Let $M$ be a complete Riemanninan three-manifold with $R_{0}\geq 0$. Let $\Su$ be a complete and stable minimal surface  diffeomorphic to the punctured disc $S^{1}\times \field{R}$. Then $\Su$ is flat and totally geodesic.
\end{Lemma}

\n For the proof we will need the next Proposition. 

\begin{Proposition}\label{P6} Assume that $R_{0}\geq 0$. Let $\Su$ be a stable surface and let $\ell_{2}$ be a smooth embedded loop dividing $\Su$ into two connected components $U_{1}$ and $U_{2}$. Let $L_{1}$, $L$ and $L_{2}$ be such that $L_{1}+L\leq Rad(\ell_{2},U_{1})$, $L_{2}\leq Rad(\ell_{2},U_{2})$ and that $L$ is less than the distance from $\ell_{2}$ to the locus of the distance function to $\ell_{2}$. Let $\ell_{1}=\partial (B(\ell_{2},L)\cap U_{1})$. Denote $l_{1}=length(\ell_{1})$ and $l_{2}=length(\ell_{2})$. Then if $l_{1}\leq l_{2}$ we have

\be\label{ine1}
(1-\frac{l_{1}}{l_{2}})^{2}\leq 4(\frac{L}{L_{1}}+\frac{L}{L_{2}}),
\ee

\n while if $l_{2}\leq l_{1}$ we have

\be\label{ine2}
(1-\frac{l_{2}}{l_{1}})^{2}\leq 4(\frac{L}{L_{1}}+\frac{L}{L_{2}}).
\ee

\end{Proposition}

\vspace{0.2cm}
\n {\bf Proof:} (of Proposition \ref{P6})

Let us give names to the regions separated by the loops $\ell_{1}$ and $\ell_{2}$. Let ${\mathcal{R}}_{2}=U_{2}$, let ${\mathcal{R}}$ be the region between $\ell_{1}$ and $\ell_{2}$ and let ${\mathcal{R}}_{1}$ be the third region, namely ${\mathcal{R}}_{1}=\Su\setminus ({\mathcal{R}}_{2}\cup {\mathcal{R}})$. We are going to consider a trial function $f$ in the operator (\ref{svo}) that we describe next. 

\vs
\begin{equation*}
\begin{array}{ll}
f(p)=f_{2}(1-\frac{dist(p,\ell_{2})}{L_{2}})&\ {\rm if\ } p\in {\mathcal{R}}_{2}\ {\rm and}\ dist(p,\ell_{2})\leq L_{2},\\
f(p)=0&\ {\rm if\ } p\in {\mathcal{R}}_{2}\ {\rm and}\ dist(p,\ell_{2})> L_{2},\\
&\\
f(p)=(f_{1}-f_{2})\frac{dist(p,\ell_{2})}{L}+f_{2}&\ {\rm if\ } p\in {\mathcal{R}},\\
&\\
f(p)=f_{1}(1-\frac{dist(p,\ell_{1})}{L_{1}})&\ {\rm if\ } p\in {\mathcal{R}_{1}}\ {\rm and}\ dist(p,\ell_{1})\leq L_{1},\\
f(p)=0&\ {\rm if\ } p\in {\mathcal{R}_{1}}\ {\rm and}\ dist(p,\ell_{1})> L_{1}.\\
\end{array}
\end{equation*}

\vs
\n Above $f_{1}$ and $f_{2}$ are constants to be given later. We use now Theorem \ref{T1} to the linear functions that define $f$ on ${\mathcal{R}}_{1}$ and ${\mathcal{R}}_{2}$. We get
\ben
\int_{{\mathcal{R}}_{1}\cup{\mathcal{R}}_{2}}|\nabla f|^{2}+\kappa f^{2}dA\leq f_{2}^{2}(2\frac{l_{2}}{L_{2}}+l_{2}'-\frac{A_{2}}{L_{2}^{2}}) + f_{1}^{2}(2\frac{l_{1}}{L_{1}}+l_{1}'-\frac{A_{1}}{L_{1}^{2}}),
\een

\n where, following Theorem \ref{T1}, $l'_{1}$ and $l'_{2}$ are the derivatives of the lengths of $\ell_{1}$ and $\ell_{2}$ when they are displaced in the inward direction to ${\mathcal{R}}_{1}$ and ${\mathcal{R}}_{2}$ respectively. The contribution from the region ${\mathcal{R}}$ can be calculated explicitly (see the proof of Theorem \ref{T1} and note that there is no locus) and gives
\ben
\int_{{\mathcal{R}}}|\nabla f|^{2}+\kappa f^{2}dA = -f'^{2}A+2f_{1}f'l_{1}-2 f_{2}f' l_{2}-f_{2}^{2}l_{2}'-f_{1}^{2}l_{1}',
\een

\n where $f'$ is the derivative of $f$ as a function of $dist(p,\ell_{2})$ and therefore equal to $f'=(f_{1}-f_{2})/L$. Adding both contributions up, discarding the negative terms containing the area factors $A_{1},A$ and $A_{2}$, and canceling the terms containing the factors $l'_{1}$ and $l'_{2}$ we obtain
\ben
\int_{\Su}|\nabla f|^{2}+\kappa f^{2}dA\leq 2f_{2}^{2}\frac{\ell_{2}}{L_{2}}+2\frac{\ell_{1}}{L_{1}}+
2\frac{(f_{1}-f_{2})}{L}(f_{1}l_{1}-f_{2}l_{2}).
\een

\n Choose $f_{1}=1$ and $f_{2}=(l_{1}+l_{2})/(2l_{2})$. We deduce
\ben
\frac{f_{1}-f_{2}}{L}=\frac{l_{2}-l_{1}}{2L l_{2}},
\een
\ben
f_{1}l_{1}-f_{2}l_{2}=\frac{l_{1}-l_{2}}{2}.
\een

\n Thus we get ($R_{0}=0$)
\ben
0\leq \int_{\Su}|\nabla f|^{2}+\kappa f^{2} dA\leq 2\frac{l_{1}}{L_{1}}+\frac{1}{2}\frac{(l_{1}+l_{2})^{2}}{l_{2} L_{2}}-\frac{(l_{1}-l_{2})^{2}}{2L l_{2}}.
\een

\n Therefore
\ben
\frac{(l_{1}-l_{2})^{2}}{2l_{2}L}\leq 2\frac{l_{1}}{L_{1}}+\frac{(l_{1}+l_{2})^{2}}{2l_{2}L_{2}}.
\een

\n As $l_{1}\leq l_{2}$ we get
\ben
\frac{(l_{1}-l_{2})^{2}}{2l_{2}L}\leq 2\frac{l_{2}}{L_{1}}+2\frac{l_{2}}{L_{2}}.
\een

\n which after multiplying by $L$ and dividing by $l_{2}$ gives the inequality (\ref{ine1}). The case when $l_{2}\leq l_{1}$ proceeds along the same lines. 

\ep 

\n {\bf Proof:} (of Lemma \ref{L2}) 

Pick a loop $\ell$ isotopic to the $S^{1}$. Observe that if  the surface $\Su$ is diffeomorphic to $S^{1} \times \field{R}$ and complete then one can take $L_{1}=L_{2}=\infty$ (independently on $\ell$ and $L$). It follows from Proposition \ref{P6} that $l_{1}=l_{2}$ for any $\ell$ and $L$ less than the distance from $\ell$ to the locus of the distance function to $\ell$. Thus for any embedded loop $\ell$ isotopic to $S^{1}$ we have
$l'=0$ where $l'$ is the derivative of $l(\ell)$ when we displace $\ell$ in its normal direction. Now if at a point $p$ in $\Su$ it is $\kappa\neq 0$ (where $\kappa$ is the Gauss curvature) then pick any two smooth loops (each one isotopic to the factor $S^{1}$ of $\Su$) such that $\ell_{1}$ passes through $p$ and $\ell_{2}$ is disjoint from $\ell_{1}$ only inside a ball around $p$ where $\kappa$ doesn't change sign. It follows by Gauss-Bonet that the integral of $\kappa$ inside the region enclosed by $\ell_{1}$ and $\ell_{2}$ must be zero which is a contradiction.\ep  

\vs
{\center \subsection{Integral curvature and isoperimetric collapse in dimension three.}\label{S3.4}}
\vs

Lemma \ref{L2} has the following interesting consequence in Riemannian geometry (the result follows indeed from \cite{Fischer-Colbrie-Schoen}). 

\vs
\begin{Proposition}\label{PT} Let $M$ be a Riemannian three-manifold with $R\geq 0$ and diffeomorphic to $S^{1}\times S^{1}\times \field{R}$. Suppose that the three-metric $g$ on $M$ is asymptotic to $h_{L}+ dx^{2}$ on one end of $\Su$ and where $h_{L}$ is a flat metric
on $S^{1}\times S^{1}$ and, asymptotic to $h_{R}+ dx^{2}$ on the other end of $\Su$. Then $h_{L}=h_{R}$.
\end{Proposition}

\n {\bf Proof:} (Sketch)

To see this consider two stable and complete tubes, isotopic to $S^{1}\times \field{R}$  for each one of the factors $S^{1}$ in $S^{1}\times S^{1}$. Such tubes must be flat and totally geodesic. But if two totally geodesic surfaces intersect along a line they must preserve the angle of intersection. It follows that $h_{R}=h_{L}$.\ep 

\begin{figure}[h]
\centering
\includegraphics[width=7cm,height=4cm]{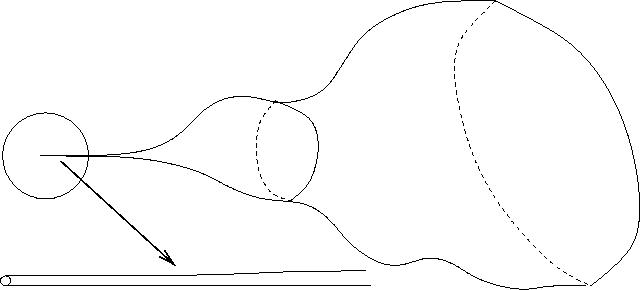}
\caption{A surface with a geometric configuration that is forbiden for stable surfaces.}
\label{Fig3}
\end{figure}

\vs
Infinite, stable and complete tubes (homeomorphic to $S^{1}\times (-1,1)$ and consequently of infinite diameter) on three-manifolds of non-negative scalar curvature must be flat, therefore (and trivially) the size of geodesic loops (isotopic to the factor $S^{1}$) keep their length constant when we translate them from one end of the tube to the other. There is a relative property to this, for stable compact tubes ($\Su$ homeomorphic to $S^{1}\times [-1,1]$), on three-manifolds of non-negative scalar curvature. The property can be interpreted as saying that: under a uniform lower bound on the area and a uniform upper bound on the diameter, then no one of the ends of $\Su$ can become arbitrarily thin (see the precise statement of the proposition below). Figure \ref{Fig3} shows a surface with a geometric configuration that is forbiden for stable surfaces. This property is central, as we will discuss briefly later in this section, to study the geometry of isoperimetric collapse on three-manifolds with uniform bounds on the $L^{p}_{g}$-norm of the Ricci curvature ($p>3/2$) and scalar curvature uniformly bounded below $(R\geq R_{0})$.  

Consider a Riemannian surface $(\Su,h)$ homeomorphic to the tube $[-1,1]\times S^{1}$. Suppose that the stability inequality 
\be\label{ESI}
\int_{\Su}|\nabla f|^{2}+\kappa f^{2}dA\geq 0,
\ee

\n holds for any $f$ vanishing at the boundary of $\Su$. Consider now a smooth loop $\ell$ embedded in $\Su$ and isotopic to any one of the two loops $\ell_{1}$ and $\ell_{2}$ that form the boundary of $\Su$. Let $L_{1}$ and $L_{2}$ be the distances from $\ell$ to $\ell_{1}$ and $\ell_{2}$ respectively and let $A_{1}$ be the area of the set of points at a distance less or equal than $L_{1}$ form $\ell$ in the component that contains $\ell_{1}$ (and similarly for $A_{2}$). Finally let $l=length(\ell)$, $l_{1}=length(\ell_{1})$ and $l_{2}=length(\ell_{2})$. It follows from Theorem \ref{T1}
that the following {\it size relation} holds
\be\label{CE}
2l(\frac{1}{L_{1}}+\frac{1}{L_{2}})\geq \frac{A_{1}}{L_{1}^{2}}+\frac{A_{2}}{L_{2}^{2}}.
\ee

\n Note that the expression is scale invariant. A direct consequence of this geometric relation is the following proposition. We will use the notation above in the statement as well as in the proof. 

\begin{Proposition}\label{NCFD} Suppose $A_{0}$ and $L_{0}$ are constants greater than zero. Then there are no sequences of Riemannian surfaces $\{(\Su_{i},h_{i})\}$  diffeomorphic to $S^{1}\times [-1,1]$, satisfying the stability inequality (\ref{ESI}) and also
\begin{enumerate}
\item $l_{1,i}\rightarrow 0$,
\item $dist(\ell_{1,i},\ell_{2,i})\leq L_{0}$,
\item $A(B(\ell_{1,i},dist(\ell_{1,i},\ell_{2,i})))\geq A_{0}>0$,
\item The pointed scaled spaces $(\Su_{i},\frac{1}{l_{1,i}^{2}}h_{i},p_{i})$, where $p_{i}\in \ell_{1,i}$, converge (in $C^{1,\alpha}$) to the flat tube $[0,\infty)\times S^{1}$.
\end{enumerate} 
\end{Proposition}

\n {\bf Proof:}

Denote the metric $(1/l_{1,i}^{2})h_{i}$ as $\tilde{h}_{i}$. Take an increasing and diverging sequence $\{d_{i}\}$, increasing slow enough in such a way that the region $B_{\tilde{h}_{i}}(\ell_{1,i},d_{i})$ (provided with the metric $\tilde{h}_{i}$) is closer and closer to the flat tube $[0,d_{i}]\times S^{1}$. Now chose as intermediate loop $\ell_{i}$ (in the setup of equation (\ref{CE})), the one at a distance $d_{i}/2$ from $\ell_{1,i}$. It follows from the scale invariance of the expression (\ref{CE}) that $l_{i}(1/L_{1,i}+1/L_{2,i})\rightarrow 0$. Thus we have
\ben
0=\lim \frac{A_{1,i}}{L_{1,i}^{2}}+\frac{A_{2,i}}{L_{2,i}^{2}}\geq \frac{A_{1,i}+A_{2,i}}{L_{0}^{2}}=\frac{A(B(\ell_{1,i},dist(\ell_{1,i},\ell_{2,i})))}{L_{0}^{2}}\geq \frac{A_{0}}{L_{0}^{2}},
\een

\n which is absurd.\ep

\vs
\n To understand the applicability of the previous proposition we will consider the following family (indexed in $\epsilon$) of Riemannian three-manifolds, $(g_{\epsilon}=dr^{2}+r^{2} d\theta^{2}+(\epsilon+r^{2})^{4}d\theta_{1}^{2},D^{2}\times S^{1})$ where $D^{2}$ is the two-dimensional unit disc and $(r,\theta)$ are polar coordinates on $D^{2}$. As $\epsilon\rightarrow 0$, it is easy to see that the (three-diemsional) volume radius at the central circle $\{r=0\}$ collapses (tends) to zero, while the (three-dimensional) volume radius at any point in the torus $\{r=1\}$ remains uniformly bounded below. Also, the isoperimetric constant tends to zero. It can be easily seen too that the integral curvature $\|Ric\|_{L^{2}_{g}}$ remains uniformly bounded above. This example (essentially due to D.Yang) clearly shows that there can be isoperimetric collapse at a finite distance from a non-collapsed region, while keeping the $L^{2}_{g}$- norm  of the Ricci curvature uniformly bounded above. Situations like this represent a serious technical inconvenience for the theory of convergence/collapse of Riemannian three-manifolds under $L^{p}_{g}$-curvature bounds \cite{Anderson}. However this inconvenience is prevented if the scalar curvature is assumed to be non-negative (situations like the previous example can be ruled with just any a priori lower bound on the scalar curvature). To see this, consider the stable tubes on $D^{2}\times S^{1}$ given by $\Su=\{(r,\theta,\theta_{1}),r\in[0,1],\theta=0,\theta_{1}\in [0,2\pi)\}$ and apply Proposition \ref{NCFD}. Indeed a direct computation shows that the scalar curvature at the central fiber ($\{r=0\}$) tends to negative infinity. A discussion of the applicability of Proposition \ref{NCFD} to three-manifolds with scalar curvature a priori bounded below and to General Relativity is given in \cite{Reiris}. 

\vs
{\center \section{A proof of Theorem \ref{T1}.}\label{ST1}}
\vs

In this section we are going to give a proof of Theorem \ref{T1}. For clarity, we state it again below.

\vs
\begin{Theorem} Let $\Su$ be a compact surface with smooth boundary consisting of a finite set of closed curves $\partial \Su=\{\ell_{1}\cup\ldots\cup \ell_{n}\}$. Consider the distance function $r(p)=dist(p,\ell_{i_{1}}\cup\ldots\cup \ell_{i_{j}})$ where $j $ is a fixed number between $1$ and $n$. Fix a distance $L$ less or equal than the distance between $\ell_{i_{1}}\cup\ldots\cup \ell_{i_{j}}$ and $\partial \Su\setminus (\ell_{i_{1}}\cup\ldots\cup \ell_{i_{j}})$ and consider the radial function $f(p)=1-\frac{r(p)}{L}$ if $r(p)\leq L$ and zero otherwise. Then we have
\be\label{svcII}
\int_{\Su}|\nabla f|^{2}+\kappa f^{2}dA\leq 2\frac{l}{L}+l'-\frac{A}{L^{2}},
\ee

\n were $l=length(\ell_{i_{1}})+\ldots +length(\ell_{i_{j}})$ and $l'=length(\ell_{i_{1}})'+\ldots+length(\ell_{i_{j}})'$ where $'$ denotes the derivative when the curve is displaced in its inward normal direction. Also $A=Area(B(\ell_{i_{1}}\cup\ldots\cup \ell_{i_{j}},L))$.
\end{Theorem}

We will give an explicit proof of this statement assuming that the distance function $r(p)=dist(p,\ell_{i_{1}}\cup \ldots\cup \ell_{i_{j}})$ has {\it regular cut locus}. We will give the definition of {\it regular  cut locus} below.  Before, let us recall the definition of cut locus. Denote by ${\mathcal{N}}(\ell_{i_{1}}\cup\ldots\cup \ell_{i_{j}})$ the set of tangent vectors $N(p)$ to $\Su$, based at points $p$ in $\ell_{i_{1}}\cup \ldots \cup\ell_{i_{j}}$, perpendicular $\ell_{i_{1}}\cup \ldots \cup\ell_{i_{j}}$ and inward pointing to $\Su$. Consider the exponential map $exp: {\mathcal{N}}(\ell_{i_{1}}\cup\ldots\cup \ell_{i_{j}})\rightarrow \Su$ that to a given $N(p)$ it assigns the point $exp(p,N(p))$, i.e. the end point of the geodesic with initial velocity $N(p)/|N(p)|$ and arc length $|N(p)|$. {\it A point $q$ is in the cut-locus of the exponential map iff either there are at least two geodesics segments, minimizing the distance between $\ell_{i_{1}}\cup\ldots\cup \ell_{i_{j}}$ and $q$ or there is one but the differential of the exponential map fails to be invertible}. Thus, to every point $q$ not in the locus of the exponential map there is a unique geodesic segment joining $q$ to $\ell_{i_{1}}\cup\ldots\cup\ell_{i_{j}}$, minimizing the distance between them (and naturally perpendicular to $\ell_{i_{1}}\cup\ldots\cup\ell_{i_{j}}$).  We will use the terminology {\it locus of the distance function $r(p)=dist(p,\ell_{i_{1}}\cup\ldots\cup\ell_{i_{j}})$} when we want to refer to the cut-locus of the exponential map  $exp:{\mathcal{N}}(\ell_{i_{1}}\cup\ldots\cup \ell_{i_{j}})\rightarrow \Su$. 

\begin{Definition} Let $\Su$ be a smooth compact surface with smooth boundary consisting of a finite set of closed curves $\partial \Su=\{\ell_{1},\ldots,\ell_{n}\}$. Consider the distance function $r(p)=dist(p,\ell_{i_{i}}\cup\ldots \cup \ell_{i_{j}})$, $1\leq j\leq n$.  We say that the distance function has regular cut-locus iff it consists of a finite union of embedded, compact and non-closed curves, such that if two of them intersect then they do only at their end points. 
\end{Definition}

\begin{figure}[h]
\centering
\includegraphics[width=11cm,height=7cm]{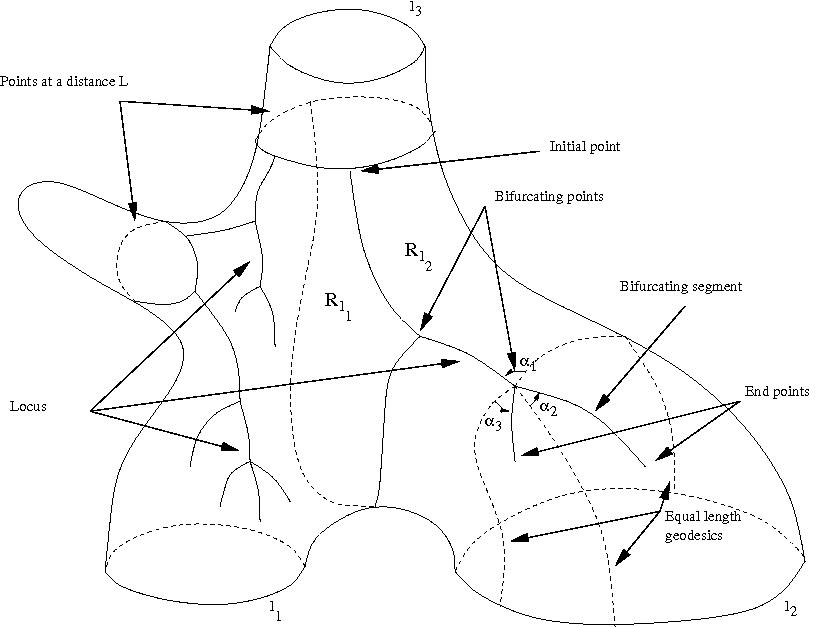}
\caption{Example of a regular cut-locus. In this case $\ell_{i_{1}}\cup\ldots\cup\ell_{i_{j}}=\ell_{1}\cup\ell_{2}$. We have signaled a bifurcating point where three bifurcating segments meet. For that point we have indicated the interior angles $\alpha_{1},\alpha_{2},\alpha_{3}$ and the three equal length geodesics meeting at it.}
\label{Fig1}
\end{figure}

\n Every embedded, compact and non-closed curve will be called a {\it bifurcating segment} and a point where two of them intersect will be called a {\it bifurcating point} (to justify this terminology see the properties they have below). The Figure \ref{Fig1} shows an example of a regular cut locus. Naturally, it may be that the cut locus has a non-regular structure but it is not the case for instance when $g$ is an analytic metric. The assumption of a regular locus does not imply any restriction on the conclusion of Theorem \ref{T1}. A regular locus has additional geometric properties that are deduced from the definition of cut-locus (as was recalled above). It is straightforwardly checked that  

\begin{enumerate}
\item  If a point $p$ belongs to the interior of a bifurcating segment then there are two and only two length minimizing geodesic segments of equal length, starting at $p$, ending at $\ell_{i_{1}}\cup\ldots\cup \ell_{i_{j}}$ and having equal interior angles with the bifurcating segment. 

\item At a bifurcating point  where $m$ bifurcating segments join, there meet $m$ (and only $m$) length minimizing geodesic segments of the same length. Two consecutive geodesic segments  enclose a region (near the bifurcating point) containing one and only one bifurcating segment and forming equal interior angles with that particular bifurcating segment. There are thus $m$ interior angles $\{\alpha_{1},\ldots,\alpha_{m}\}$ formed by these geodesics and the $m$ bifurcating segments, whose sum, $\sum_{i=1}^{i=m}\alpha_{i}$, is equal to $\pi$. 
\end{enumerate}

\n An illustration of these two properties on a regular locus is given in Figure \ref{Fig1}. 

To obtain the upper bound (\ref{svcII}) we will proceed as follows. We divide first $B(\ell_{i_{1}}\cup\ldots\cup\ell_{i_{j}},L)$ (which is the support of $f$) into $j$ disjoint regions, each one homeomorphic to $S^{1}\times [-1,1)$, plus the locus of $r(p)$. We find then a general expression for the integration of the left hand side of (\ref{svcII}) applicable to each one of the regions. We add finally all the contributions and use a combinatoric argument to conclude, after some cancelations, the upper bound (\ref{svcII}). We explain now how these regions are defined. 

Consider one of the loops $\ell$ from $\ell_{i_{1}},\ldots,\ell_{i_{j}}$. Let $p$ be a point in $\ell$ and let $\gamma_{p}$ be the geodesic emanating from $p$ an normal to $\ell$. It may be that $\gamma_{p}$ reaches the locus ${\mathcal{C}}$ of the distance function $r(p)$, at a point $q$ with arc length less than $L$ or it does not. If it does then we remove from $\gamma_{p}$ all points after $q$ and if not we consider the point $q$ with arc length equal to $L$ and remove from $\gamma_{p}$ all the points after $q$. Thus to every $p$ there is a geodesic segment, denoted by $\tilde{\gamma}_{p}$, emanating from $p$, perpendicular to $\ell$, and that terminates either at the locus of $r(p)$ or at a point with arc length equal to $L$. Naturally, the map $p\rightarrow q$, where $q$ is the end point of $\tilde{\gamma}_{p}$ is continuous and the set ${\mathcal{R}}_{\ell}=\cup_{p\in \ell} (\tilde{\gamma}_{p}\setminus \{q\})$ is homeomorphic to $S^{1}\times [-1,1)$. The regions ${\mathcal{R}}_{\ell},\ \ell\in \{\ell_{i_{1}},\ldots,\ell_{i_{j}}\}$ and the locus ${\mathcal{C}}$ of $r(p)$ form the desired partition of $B(\ell_{i_{1}}\cup\ldots\cup\ell_{i_{j}},L)$.   

At everyone of the regions ${\mathcal{R}}_{\ell}$ we will consider geodesic normal coordinates which are defined as is usual as follows. Let $p_{\ell}$ be a fixed point in $\ell$. Let $p\in \ell$ be any other point. The orientation of $\Su$ and the inner normal direction $e_{2}(p)$ to $\ell$ at $p$ define a direction $e_{1}(p)$ on $\ell$ at $p$ by imposing that the orthonormal pair $\{e_{1},e_{2}\}$ defines the same orientation as $\Su$. Given $p$ in $\ell$, define $\theta(p)$ as the length (in $\ell$) of the segment joining  $p_{\ell}$ to $p$ (following the direction $e_{1}$ before). Given a point $q$ in ${\mathcal{R}}_{\ell}$ consider the unique geodesic segment $\tilde{\gamma}_{p(q)}$ (see the notation above) passing through $q$. Let $r(q)$ be the distance between $p(q)$ and $q$ along $\tilde{\gamma}_{p(q)}$. Now to every point $q$ in ${\mathcal{R}}_{\ell}$ we associate unique coordinates $(r,\theta)=(r(q),\theta(p(q)))$. In these coordinates the metric is written as $h=dr^{2}+\phi^{2}d\theta^{2}$, the Gaussian curvature is given by $\phi''=-\kappa \phi$ and the element of area by $dA=\phi dr d\theta$.  

We move now forward to obtain a general formula for the integral on the left hand side of (\ref{svcII}) on each one of the regions ${\mathcal{R}}_{\ell_{i_{k}}}$, $k=1,\ldots,j$. Consider again a component $\ell$ of $\ell_{i_{1}}\cup \ldots \cup\ell_{i_{j}}$. We will restrict the discussion to the region ${\mathcal{R}}_{\ell}$. We will divide further ${\mathcal{R}}_{\ell}$ into two regions denoted by ${\mathcal{R}}_{\ell}^{I}$ and ${\mathcal{R}}_{\ell}^{II}$. ${\mathcal{R}}_{\ell}^{I}$ is defined as the union of the set of geodesic segments $\tilde{\gamma}_{p}\setminus \{q(p)\}$ of length $L$ and ${\mathcal{R}}_{\ell}^{II}$ is defined as the union of the set of geodesic segments $\tilde{\gamma}_{p}\setminus \{q(p)\}$ whose end point lies in the locus ${\mathcal{C}}$. 

\vs
\n {\it Integration on the $\Reg^{I}$ region.}

\vs 
We note first that $\Reg^{I}$ naturally consists of a disjoint set of connected subregions $\Reg^{I,k}$, $k=1,\ldots m(\ell)$, each homeomorphic to $[-1,1]\times [-1,1)$, or, in terms of the coordinates $(r,\theta)$, subregions of the form $[\theta_{i},\theta_{f}]\times [0,L)$. Below, $\Reg^{I,k}$ will mean one of such connected components.  We want to compute
\be\label{svfc}
\int_{\Reg^{I,k}}|\nabla f|^{2}+\kappa f^{2}dA=\int_{\theta_{i}}^{\theta_{f}}\int_{0}^{L}f'^{2}\phi-f^{2}\phi'' dr d\theta.
\ee

\n where above $'$ is the derivative with respect to $r$ and $f$, as in the statement of the Theorem, is $f(p)=1-r(p)/L$. We will integrate by parts in the variable $r$ and then integrate in the variable $\theta$. For convenience we will perform the integration by parts with an upper limit of integration of $r(\theta)$ instead of $L$ (where $r(\theta)$ is some positive function). Noting that $f-=-1/L$ we compute
\be\label{IBP}
-\int_{0}^{r(\theta)}\phi''f^{2}dr=-\phi'f^{2}|_{0}^{r(\theta)} +2\phi f f'|_{0}^{r(\theta)}-2f'^{2}\int_{0}^{r(\theta)}\phi dr.
\ee

\n We now make $r(\theta)=L$ and integrate the expression in $\theta$. Recalling $f(L,\theta))=0$ and $f(0,\theta)=1$ we get for the first term of the right hand side of (\ref{IBP})
\ben
\int_{\theta_{i}}^{\theta_{f}}-\phi'f^{2}|_{0}^{L}d\theta=\int_{\theta_{i}}^{\theta_{f}}\phi'(0,\theta)d\theta=l_{k}'.
\een

\n For the second term instead we get
\ben
\int_{\theta_{i}}^{\theta_{f}}2\phi f f'|_{0}^{r(\theta)}d\theta=2\frac{l_{k}}{L},
\een

\n and finally for the last term we obtain
\ben
-2f'^{2}\int_{\theta_{i}}^{\theta_{f}}\int_{0}^{r(\theta)}\phi dr d\theta=-2\frac{A_{k}}{L^{2}}.
\een

\n The last three terms together amount for the second term on the right hand side of equation (\ref{IBP}). The first term is equal to $A_{k}/L^{2}$ and thus we get the final result
\be\label{svfc}
\int_{\Reg^{I,k}}|\nabla f|^{2}+\kappa f^{2}dA=2\frac{l_{k}}{L}+l_{k}'-\frac{A_{k}}{L^{2}}.
\ee

\n Adding up the contributions from all the subregions $\Reg^{I,k}$, $k=1,\ldots,m(\ell)$ gives
\be\label{svfcII}
\int_{\Reg^{I}}|\nabla f|^{2}+\kappa f^{2}dA=2\frac{l}{L}+l'-\frac{A}{L^{2}},
\ee

\n where $l=\sum_{k=1}^{k=m(\ell)} l_{k}$ is the length of the base of the region $\Reg^{I}$ (in $\ell$), $l'=\sum_{k=1}^{k=m(\ell)}l'_{k}$ is the derivative of the length of the base when it is displaced in the normal inward direction and $A=\sum_{k=1}^{k=m(\ell)} A_{k}$ is the total area of $\Reg^{I}$.

\vs
\n {\it Integration on the $\Reg^{II}$ region.}

\vs We note first that any geodesic $\tilde{\gamma}_{p}$ belonging to the region $\Reg^{II}$ terminates in the locus ${\mathcal{C}}$ and that ${\mathcal{C}}$ consists of a finite set of bifurcating segments $\{\Gamma_{1},\ldots,\Gamma_{m(\ell)}\}$ (in that region). Therefore if we fix a side of a bifurcating segment $\Gamma_{k}$ we can consider the subregion $\Reg^{II,k}$ formed by all geodesic segments $\tilde{\gamma}_{p}\setminus \{q(p)\}$ ending at $\Gamma_{k}$ from the chosen side. Each $\Reg^{II,k}$ is naturally homeomorphic to $[-1,1]\times [-1,1]$ or, in terms of the coordinates $(r,\theta)$, to a domain of the form $\{(r,\theta)/ \theta\in [\theta_{i},\theta_{f}], 0\leq r\leq r(\theta)\}$ where $r$ is a differentiable and positive real function. The union $\cup_{k=1}^{k=m(\ell)} \Reg^{II,k}$ is all the region $\Reg^{II}$ and the interiors of the regions $\Reg^{II,k}$ are disjoint. Note that every bifurcating segment has two subregions associated to it (one for each one of the sides).

Consider a bifurcating segment $\Gamma_{k}\in \{\Gamma_{1},\ldots,\Gamma_{m(l)}\}$ and $\tilde{\gamma}_{p}$ a geodesic segment ending at $q\in \Gamma_{k}$ on it, that we can think is inside the subregion $\Reg^{II,k}$. We will define now the {\it interior angle} formed by $\Gamma_{k}$ and $\tilde{\gamma}_{p}$ at $q$. For this consider the vector $e_{1}(q)$ normal to $\Gamma_{k}$ at $q$ and pointing outward to $\Reg^{II,k}$. Then there is a unique unit vector $e_{2}(q)$ tangent to $\Gamma_{k}$ at $q$ such $\{e_{1}(q),e_{2}(q)\}$ defines the same orientation as the one of $\Su$. The interior angle $\alpha$ between $\Gamma_{k}$ and $\tilde{\gamma}_{p}$ at $q$ is defined as the one (between $0$ and $\pi$) formed by the velocity vector $\tilde{\gamma}_{p}'$ of $\tilde{\gamma}_{p}$ at $q$ and $e_{2}(q)$.       

We consider now a subregion $\Reg^{II,k}=\{(r,\theta)/ \theta\in [\theta_{i},\theta_{f}], 0\leq r\leq r(\theta)\}$ and perform the integration on it. Following the same calculation as in the subregions $\Reg^{I}$ (see equation (\ref{svfc}) and (\ref{IBP})) we get
\be\label{svcrII}
\int_{\Reg^{II,k}} |\nabla f|^{2}+\kappa f^{2}dA=2\frac{l_{k}}{L}+l_{k}'-\frac{A_{k}}{L^{2}} -\int_{\theta_{i}}^{\theta_{f}}2\frac{f(r(\theta))\phi(r(\theta))}{L}+\phi'(r(\theta))f^{2}(r(\theta))d\theta.
\ee

\n Let $q_{i(f)}$ be the points $(r(\theta_{i(f)}),\theta_{i(f)})$. Given $q\in \Gamma_{k}$ let $s(q)$ be the distance from $q$ to $q_{i}$ along $\Gamma_{k}$. In performing the integral (\ref{svcrII}) we will use the identity
\be\label{chvar}
\frac{d \alpha}{d\theta}=\phi'+k_{in}\frac{ds}{d\theta},
\ee

\n where $k_{in}$ is the (inward) mean curvature of $\Gamma_{k}$ in $\Su$.  Let us postpone the proof of this identity until later and use it in equation (\ref{svcrII}). The last term in (\ref{svcrII}) splits into two terms. We get first
\be\label{eq0}
-\int_{\theta_{i}}^{\theta_{f}}\phi'f^{2}d\theta=-\int_{\theta_{i}}^{\theta_{f}}\frac{d\alpha}{d\theta}f^{2}+\int_{s_{i}(=0)}^{s_{f}}k_{in}ds=
\ee
\be\label{eq}
=-\alpha f^{2}|_{\theta_{i}}^{\theta_{f}}+\int_{\theta_{i}}^{\theta_{f}}2\alpha f\frac{df}{d\theta}d\theta +\int_{s_{i}}^{s_{f}}k_{in}ds.
\ee

\n Note that 
\ben
\frac{df}{d\theta}d\theta=\frac{df}{d r}\frac{dr}{ds}ds,
\een

\n also that $df/dr=-1/L$ and (is easy to see) $dr/ds=-\cos \alpha$. Putting this in equation (\ref{eq}) gives  
\be\label{eq1}
-\int_{\theta_{i}}^{\theta_{f}}\phi'f^{2}d\theta=-\alpha_{f} f^{2}(q_{f})+\alpha_{i}f^{2}(q_{i})+\frac{2}{L}\int_{s_{i}}^{s_{f}}\alpha \cos \alpha f ds+\int_{s_{i}}^{s_{f}}k_{in}d s.
\ee

\n We compute now the term
\ben
-\frac{2}{L}\int_{\theta_{i}}^{\theta_{f}}f\phi d\theta,
\een

\n in the equation (\ref{svcrII}). It is simple to see that $\phi d\theta=\sin \alpha ds$ and putting this in the previous equation gives
\be\label{eq2}
-\frac{2}{L}\int_{s_{i}}^{s_{f}}f\sin\alpha ds.
\ee

\n Using the calculations (\ref{eq1}) and (\ref{eq2})  in equation (\ref{svcrII}) gives the final result
\be\label{finalr}
\int_{\Reg^{II,k}}|\nabla f|^{2}+\kappa f^{2}dA=
\ee
\ben
=2\frac{l_{k}}{L}+l_{k}'-\frac{A_{k}}{L^{2}}+\alpha_{i}f^{2}(q_{i})-\alpha_{f}f^{2}(q_{f})-\frac{2}{L}\int_{s_{i}}^{s_{f}}f(\sin \alpha -\alpha\cos \alpha) ds +\int_{s_{i}}^{s_{f}}k_{in}ds.
\een

Let us discuss now how the contributions from all the subregions $\Reg^{II,k}$ add up to give, together with (\ref{svfcII}), the inequality (\ref{svcII}). Given a bifurcating segment $\Gamma_{k}$ consider the two regions $\Reg^{II,k}$ and $\BReg^{II,k}$ that face $\Gamma_{k}$ (note that the loop $\bar{\ell}$ can be different than $\ell$). We are going to consider 
\ben
\int_{\Reg^{II,k}\cup \BReg^{II,k}}|\nabla f|^{2}+\kappa f^{2}dA,
\een

\n using equation (\ref{finalr}). The variables $r,\theta,s,\alpha$ for the region $\BReg^{II,k}$ will be denoted with a bar. We note that at a point $q\in \Gamma_{k}$ it is
$\bar{\alpha}(q)=\pi-\alpha(q)$ and $\alpha_{i}=\bar{\alpha}_{f}$, $\alpha_{f}=\bar{\alpha}_{i}$. Therefore to meaningfully compare the expression (\ref{finalr}) for the regions $\Reg^{II,k}$ and $\BReg^{II,k}$ we need to use in one of them the interior angle and in the other the complementary. We will see now how the expression 
\be\label{CCC1}
\bar{\alpha}_{i}f^{2}(\bar{q}_{i})-\bar{\alpha}_{f}f^{2}(\bar{q}_{f})-\frac{2}{L}\int_{\bar{s}_{i}}^{\bar{s}_{f}}f(\sin \bar{\alpha} -\bar{\alpha}\cos \bar{\alpha}) d\bar{s} +\int_{\bar{s}_{i}}^{\bar{s}_{f}}\bar{k}_{in}d\bar{s},
\ee

\n changes when instead of using $\bar{\alpha}(\bar{\theta})$ we use $\alpha(\bar{\theta})$. We have $d\alpha/d\bar{\theta}=-d\bar{\alpha}/d\bar{\theta}$, $cos \bar{\alpha}=-\cos \alpha$ and $\sin \bar{\alpha}=\sin \alpha$. Thus repeating the calculation from (\ref{chvar}) to (\ref{eq2}) we get for (\ref{CCC1}) the expression
\ben
\alpha_{i}f^{2}(q_{i})-\alpha_{f}f^{2}(q_{f})-\frac{2}{L}\int_{\bar{s}_{i}}^{\bar{s}_{f}}f(\sin \alpha -\alpha\cos \alpha) d\bar{s} +\int_{\bar{s}_{i}}^{\bar{s}_{f}}\bar{k}_{in}d\bar{s}.
\een

\n Thus we get 
\be\label{finalrb}
\int_{\BReg^{II,k}}|\nabla f|^{2}+\kappa f^{2}dA=
\ee
\ben
=2\frac{\bar{l}_{k}}{L}+\bar{l}_{k}'-\frac{\bar{A}_{k}}{L^{2}}+\alpha_{i}f^{2}(q_{i})-\alpha_{f}f^{2}(q_{f})-\frac{2}{L}\int_{\bar{s}_{i}}^{\bar{s}_{f}}f(\sin \alpha -\alpha\cos \alpha) d\bar{s} +\int_{\bar{s}_{i}}^{\bar{s}_{f}}\bar{k}_{in}d\bar{s}.
\een

\n Noting that a point $q$ in $\Gamma_{k}$ it is $k_{in}(q)=-\bar{k}_{in}$, we deduce then that when adding the contributions (\ref{finalr}) and (\ref{finalrb}) the last terms cancel out. We get thus the result 
\be\label{finalrf}
\int_{\Reg^{II,k}\cup \BReg^{II,k}}|\nabla f|^{2}+\kappa f^{2}dA=
\ee
\ben
=2\frac{l_{k}}{L}+l_{k}'-\frac{A_{k}}{L^{2}}+2\frac{\bar{l}_{k}}{L}+\bar{l}_{k}'-\frac{\bar{A}_{k}}{L^{2}}+2\alpha_{i}f^{2}(q_{i})-2\alpha_{f}f^{2}(q_{f})-\frac{4}{L}\int_{s_{i}}^{s_{f}}f(\sin \alpha -\alpha\cos \alpha) ds.
\een

\n Various elements we have to remark on this equation. First the expression $\sin \alpha -\alpha \cos \alpha$ is greater or equal than zero (for $0\leq \alpha\leq \pi$), and therefore the last term on the right hand side of equation (\ref{finalrf}) is always non-positive. Secondly, the term $2\alpha_{i}f^{2}(q_{i})-2\alpha_{f}f^{2}(q_{f})$ depends, naturally, on the choice of which one, of the two ends of $\Gamma_{k}$, is the initial ($q_{i}$) and which one the final ($q_{f}$). In other words this term depends on an a priori selection of a direction or orientation in $\Gamma_{k}$. From (\ref{svfcII}) and (\ref{finalrf}) we see that inequality (\ref{svcII}) would be satisfied as long as, for any regular locus, it is always possible to make a choice of a direction at each bifurcating segment in such a way that the sum of the terms $2\alpha_{i}f^{2}(q_{i})-2\alpha_{f}f^{2}(q_{f})$ is non-positive. As it turns out this is possible to do, we start the discussion next.     

Let $\{\Gamma\}$ be the set of bifurcating points of the locus ${\mathcal{C}}$. Suppose that to every $\Gamma$ it was prescribed a direction or orientation $\{\vec{\Gamma}\}$. This done, every oriented bifurcating segment $\vec{\Gamma}$ will have an initial point $q_{i}(\vec{\Gamma})$ and a final pint $q_{f}(\vec{\Gamma})$. Similarly will have initial angles $\alpha_{i}(\vec{\Gamma})$ and final angles $\alpha_{f}(\vec{\Gamma})$. We would like to discuss how the orientations or directions of the bifurcating segments have to be prescribed to have 
\be\label{sum}
\sum_{\{\vec{\Gamma}\}}2\alpha_{i}(\vec{\Gamma})f^{2}(q_{i}(\vec{\Gamma}))-2\alpha_{f}(\vec{\Gamma})f^{2}(q_{f}(\vec{\Gamma}))\leq 0.
\ee

There are three kinds of boundary points of bifurcating segments. First are those, as we introduced before, that are shared by at least two bifurcating segments. We have called them {\it bifurcating points}. Secondly there are those which are the boundary point of a single bifurcating segment. We can distinguish two kinds of these: those that lie at a distance $L$ from $\ell_{i_{1}}\cup\ldots\cup \ell_{i_{j}}$ and those that lie at a distance less than $L$. Points of the first kind will be called {\it initial points} and points of the second kind will be called {\it end points}. We can think that the locus ``begins" at {\it initial points} and ``ends" at {\it end points} (see Figure \ref{Fig1}). 

Consider a bifurcating segment $\Gamma_{k}$ having an end-point $q$ as one of its boundary points. If we prescribe on $\Gamma_{k}$ the direction pointing towards the end-point $q$ then the end-point would be the final point $q=q_{f}$ and we would have the contribution $-2\alpha_{f}f^{2}(q_{f})$ to the sum (\ref{sum}) which is non-positive and indeed zero as $\alpha_{f}=0$. If instead we prescribe the direction pointing outwards then the end-point will be the initial point ($q=q_{i}$) and will have the contribution $2\alpha_{i}f^{2}(q_{i})$ to the sum (\ref{sum}) which is positive and indeed equal to $2\pi f^{2}(q_{i})$ as $\alpha_{i}=\pi$. {\it We conclude thus that the end-point $q$ would contribute negatively to the sum (\ref{sum}) as long as we prescribe the direction in $\Gamma_{k}$ that points towards it}.          

Consider now a bifurcating point $q$ and suppose there are $m$ bifurcating segments 

\n $\{\Gamma_{k_{1}},\ldots,\Gamma_{k_{m}}\}$ meeting at $q$. There are thus $m$ interior angles $\{\alpha_{k_{1}},\ldots,\alpha_{k_{m}}\}$ adding $\pi$, 

\n $\sum_{i=1}^{i=k}\alpha_{k_{i}}=\pi$. Let now $\Gamma_{k_{i}}$ be one of the bifurcating segments and $\alpha_{k_{i}}$ the interior angle associated to it. We note that $q$ and $\alpha_{k_{i}}$, will be equal to the initial point ($q=q_{i}$) and angle $\alpha_{k_{i}}=\alpha_{i}$ respectively, if we prescribe on $\Gamma_{k}$ the direction pointing outward from $q$. But if we prescribe on $\Gamma_{k_{i}}$ the direction point towards $q$ then $q$ will be the final point ($q=q_{f}$) and the final angle $\alpha_{f}$ will be equal to $\alpha_{f}=\pi-\alpha_{k_{i}}$. {\it Recalling that $\sum_{i=1}^{i=m}\alpha_{k_{i}}=\pi$ we conclude thus that the bifurcating point $q$ would contribute non-positively to the sum (\ref{sum}) as long as we prescribe at least one of the directions on the bifurcating segments $\Gamma_{k_{i}}$ pointing towards $q$.} 

Consider now an initial-point $q$. We note that because at any initial point $q$ it is $f(q)=0$ the contribution from any initial point to the sum (\ref{sum}) would be equal to zero. {\it We conclude thus that if $\Gamma_{k}$ is a bifurcating segment having an initial-point $q$ as one of its boundary points then it does not matter which prescription we make on the direction on $\Gamma_{k}$ that $q$ will contribute to zero to the sum (\ref{sum}).}  

Based on the above three conclusions we give now a constructive procedure to assign directions to the bifurcating segments in such a way as to to satisfy inequality (\ref{sum}). Given $\ell$ from $\{\ell_{i_{1}},\ldots,\ell_{i_{j}}\}$ we will assign first directions to the bifurcating segments in the region $\Reg$. Bifurcating segments can be shared by two different regions, because of this, we will explain over the end how to fix a consistent orientation to all bifurcating segments.    

\begin{figure}[h]
\centering
\includegraphics[width=7cm,height=7cm]{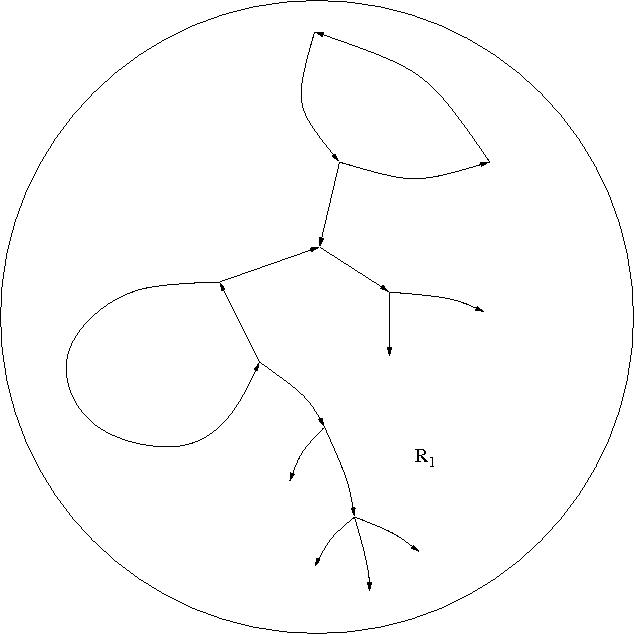}
\caption{The planar graph arising from the region ${\mathcal{R}}_{\ell_{1}}$ in Figure \ref{Fig1}. The orientation on the edges have been provided according to the procedure in Steps $1,2,3$.}
\label{Fig2}
\end{figure}

There is an important realization of the closure of $\Reg$ as a planar domain inside the unit disc. Indeed recall that $\Reg$ is homeomorphic to $S^{1}\times [-1,1)$, therefore we can find a homomorphism $\varphi_{\ell}$ (into the image) from the closure of $\Reg$into the unit disc in the plane in such a way that $\ell$ is mapped diffemorphically into the unit circle and the other boundary 
component of the closure of $\Reg$ gets mapped into {\it a planar graph such that no edge is enclosed by a closed chain of edges}. This last property turns out to be important. A representation of one of such graphs is given in Figure \ref{Fig2}.   

We will use the standard terminology ``edge" and ``vertex". A vertex which corresponds (under $\varphi_{\ell}$) to an end-point in ${\mathcal{C}}$ will be called an end-vertex and similarly, a point which corresponds (under $\varphi_{\ell}$) to an initial-point will be called an initial-vertex. The procedure to assign a direction to the edges is as follows.  

\begin{enumerate}
\item Step 1: For every end-vertex pick the only edge having it as a boundary component and assign to that edge the direction that points towards the end-vertex. Eliminate then all such edges from the graph. Of the remaining vertices we look now at those that 
\begin{enumerate}
\item are a boundary component of at least one edge that was eliminated,
\item are not an initial vertex,
\item are the boundary component of only one edge that was not deleted.
\end{enumerate}

\n We will think such vertices as new end-vertices. We repeat the procedure until we get a graph having no vertices satisfying the three items above. If such graph is empty then we are done, if not we work with it in the next Step. Note that any vertex of the resulting graph that is a boundary component of only one edge must be necessarily an initial-vertex.  
 
\item Step 2: We look now at all the closed chains of edges. Every one of such chains encloses a disk (with no edges inside). Assign to each edge the natural orientation (direction) that inherits as the boundary of a disk in the plane. Eliminate all of such closed chains. The resulting graph will have no closed chains. If the resulting graph is empty then we are done, if not we work with it in the next Step. Note that with such assignment of directions, every vertex of a closed chain has at least one direction pointing to it. Because of this and the note in the previous Step $1$, we will think that any vertex of the resulting graph which is a boundary component of only one edge is an initial-vertex.     

\item Step 3: We look now at all the vertices which are a boundary component of only one edge (and therefore are initial). If the edge of one of such points is isolated, i.e. it does not share a boundary component with any other edge, then assign to that edge any direction and eliminate it from the graph. If not, then assign to that edge the direction that points outwards from the initial-vertex. Eliminate from the graph all such edges. If a vertex of the resulting graph is a boundary component of only one edge then think it as a new initial-vertex. We repeat the procedure until there are no more edges in the graph.   
    
\end{enumerate}  

\n It is directly checked that every vertex which is not an initial vertex gets at least one direction pointing to it. 

Every component $\Reg$, $\ell\in \{\ell_{i_{1}},\ldots,\ell_{i_{j}}$ has associated a graph and some graphs can share edges. To assign a direction to every edge of every graph without incurring in inconsistencies we proceed as follows. Order the graphs as $G_{\ell_{i_{1}}},\ldots,G_{\ell_{i_{j}}}$. Pick $G_{\ell_{i_{1}}}$ and assign to every one of its edges a direction according to the procedure in Steps $1,2,3$ above. Pick $G_{\ell_{i_{2}}}$ and assign to every edge in it that is not shared with $G_{1}$ a direction that would have from applying the procedure in Steps $1,2,3$ above to the whole graph $G_{\ell_{i_{2}}}$. Pick $G_{\ell_{i_{3}}}$ and assign to every edge not shared with $G_{1}$ or $G_{2}$ the direction that would have from applying the procedure in Steps $1,2,3$ above to the whole graph $G_{\ell_{i_{3}}}$. Continue like this until $G_{\ell_{i_{n}}}$. Naturally, every vertex that is not an initial vertex gets at least one direction pointing to it.     

Let us finish the proof of the Theorem by explaining how to obtain equation (\ref{chvar}). Consider two geodesics $\tilde{\gamma}_{p_{0}}$ and $\tilde{\gamma}_{p_{1}}$ with $\theta(p_{0})=\theta_{0}$, $\theta(p_{1})=\theta_{1}$ and reaching the bifurcating segment $Gamma$ at the points $q(\theta_{0})$, $q(\theta_{1})$. Let the interior angle at $q(\theta_{0})$  be $\alpha(\theta_{0})$ and assume that $\alpha< \pi/2$ (the case $\alpha\geq \pi/2$ proceeds in similar ways). Let the length of the geodesic $\tilde{\gamma}_{p_{1}}$ be $r(\theta_{1})$. Consider the segment formed by the points in all the geodesics $\tilde{\gamma}_{p}$ with $\theta_{0}<\theta(p)<\theta_{1}$ and arc length equal to $r(\theta_{1})$. We will call such segment {\it segment I}.  Consider in addition the segment along the geodesic $\tilde{\gamma}_{p_{0}}$ formed by the points $q$ with length $r(q)$ between $r(\theta_{0})$ and $r(\theta_{1})$. We will call such segment {\it segment II}. Finally consider the segment along the bifurcating segment between the points $q(\theta_{0})$ and $q(\theta_{1})$. We will call such segment {\it segment III}. Consider the triangle formed by the three segments and denote by $\beta(\theta_{1})$ the angle formed by the first and third segments. Note that the angle formed by the first and second segments is $\pi/2$. By Gauss-Bonnet we get
\be\label{FF}
\alpha(\theta_{0})+\beta(\theta_{1})-\frac{\pi}{2}+\int_{\theta_{0}}^{\theta_{1}}\phi'(r(\theta_{1}),\theta)d\theta +\int_{s(\theta_{0})}^{s(\theta_{1})}k_{in}ds=\int_{\triangle}\kappa dA.
\ee

\n Differentiating with respect to $\theta_{1}$ at $\theta_{1}=\theta_{0}$ gives 
\ben
\frac{d \alpha}{d\theta}=-\frac{ d\beta}{d\theta}=\phi'+k_{in}\frac{ds}{d\theta},
\een

\n because the derivative of the right hand side of equation (\ref{FF}) at $\theta_{1}=\theta_{0}$ is equal to zero.\ep

\end{document}